\documentstyle[12pt]{article}
\begin{document}
\def \Bbb{\bf}
 \def \R{{\Bbb R}}
\def \P{{\Bbb P}}
\def \S{{\Bbb S}}
\def \A{{\cal A}}
\def \Vect{{\rm Vect}}
\def \V{{\cal V}}
\def \Vp{\V_P}
\def \Vq{\V_Q}
\def \D{{\cal D}}
\def\d{{\cal D}}
\def \bDp{\overline{\D_P}}
\def \bDq{\overline{\D_Q}}
\def \Dp{\D_P}
\def  \Dq{\D_Q}
\def \dX{\delta_X}
\def \E{{\cal E}}
\def \ep{\E_{P}}
\def \eq{\E_{Q}}
\def \Sp{\S_{P,Q}}
\def \C{{\cal C}}
\def \Cp{\C_P}
\def \Cq{\C_Q}
\def \mp{m_{P,Q}}
\def\arrow{\to}
\def\semi{\stackrel{.}{+}}
\def\follows{\to}
\def\Vr{\V_r}
\def\F{{\cal F}}
\def\Fun{{\rm Fun}}
\def \ad{{\rm ad}}
\def \be{\begin{equation}}
\def \ee{\end{equation}}
%%%%%%%%%%%%%%%%%%%%%%%%%%%%%%%%%%%%%%%%%%%%%%%%%%%%%%%%%%%%%%%%%%%%%%%%%%%%%%
%\textheight=20truecm
%\textwidth=13truecm

\newtheorem{Theorem}{THEOREM}
 \newtheorem{Prop}{PROPOSITION}
 \newtheorem{Proposition}{PROPOSITION}
 \newtheorem{Lemma}{LEMMA}
 \newtheorem{lemma}{LEMMA}
 \newtheorem{Corollary}{COROLLARY} 
 \title{Quantization of Poisson Structures on $\R^2$}
 \author{Dmitry Tamarkin }
 \date{}
\abstract{An 'isomorphism' between the 'moduli space' of 
star products on $\R^2$ and 'moduli space' of all 
formal Poisson structures on $\R^2$ is established. }

 \maketitle
This is a preliminary version of the paper!

The problem of quantization of Poisson structures originates from \cite{FFF}.
It is well known that any Poisson structure on a two-dimensional manifold is quantizable. In this paper 
we establish an 'isomorphism' between the 'moduli space' of 
star products on $\R^2$ and 'moduli space' of all 
formal Poisson structures on $\R^2$ by construction of a 
map from Poisson structures to star products. Certainly, this isomorphism 
follows from the Kontsevich formality conjecture \cite{Kon}. Most likely,
our map can be used as a first step in constructing 
the $L_{\infty}$ quasiisomorphism in the formality conjecture for $\R^2$.
The author would like to thank Boris Tsygan and Paul Bressler for the attention
and helpful suggestions.
   
The set of all star-products $\S$ is acted upon by the group
$\d\semi{\rm Diffeo}\R^2$, where $\d$ is the group of operators of the form
$1+hD_1+h^2D_2+\ldots$ with $D_k$ to be arbitrary differential operators.
The set of all formal Poisson structures $\P$ consists of
formal series in $h$ with bivector fields as the coefficients. Formal Poisson
structures are acted upon by the group ${\rm Diffeo\R^2}\semi{\rm exp}(h{\rm Vect}[[h]])$, where ${\rm Vect}$ is the Lie algebra of vector fields on $\R^2$.
These actions define  equivalence relations. We want to have a pair of maps
$f_1:\S\arrow\P$ and $f_2:\P\arrow\S$ such that 
$$
f_1\circ f_2(x)\equiv x\quad
f_2\circ f_1(x)\equiv x, 
$$
\be
x\equiv y\follows f_{1,2}x\equiv f_{1,2}y.\label{equiv}
\ee
 By a map from $\S$ we mean a differential expression in terms of the coefficients of the bidifferential operators
corresponding to the star products. Maps from $\P$ are defined similarly.

We can replace $\S$ by a subspace. Let $P,Q$ be a non degenerate pair of
(real) polarizations of $\R^2$ Define a subset $\S_{P,Q}$ of $\S$ in the following
way 
$m\in \Sp$ iff $m(f,g)=fg$ if $f$ is constant along $P$ or $g$ is constant along
$Q$.
\begin{Prop} 
Let $x,y$ be a nondegenerate coordinate system on $\R^2$ such that $x$
is constant along $Q$ and $y$ is constant along $P$. Then there exists a unique
map $
\S\arrow\d:m\mapsto U(m)=1+hV(m)$ such that
$$1) \mp(m)=U^{-1}(m(Uf,Ug))\in\Sp$$ 
\be2) Ux=x,\ Uy=y,\ U1=1.\label{polar}\ee
$U$ is uniquely defined by the condition
$U(x^{*m}*y^{*n})=x^my^n$ (where star denotes the star product $m$).
\end{Prop}
We denote by $\mp:\S\to\Sp$ the map which sends $m$ to $\mp(m)$. 
Further, $x,y$ will mean the same as in Proposition 1.
Thus, it is enough to find maps $p_1:\Sp\arrow\P$ and $p_2:\P\arrow\Sp$
 with the same properties as  $f_1,f_2$ have. Indeed , put
$f_2=i\circ p_2$ and $f_1=p_1\circ \mp$ (here $i: \Sp\arrow \S$ is
the injection).
 
The following theorem gives an explicit construction for $p_2$ which appears
to be a bijective map so that we can put $p_1=p_2^{-1}$. Denote by 
$\Cp$
(resp. $\Cq$) the space of functions, constant along $Q$ (respectively $P$). Denote by $\Vp$
(resp. $\Vq$) the space of
vector fields preserving the polarizations and tangent to $P$ (resp.$Q$).
Denote by $\Dp$ the subalgebra of the algebra of the differential operators
consisting of the operators $D$ such that $D(\Cq)\subset \Cq$ and
$D(fg)=fD(g)$ if $f\in \C_P$.  Denote by $\Dq$ the same algebra, where
$P$ and $Q$ are interchanged. In the coordinates $x,y$
we have $\Cp=\{f(x)\}$, $\Vp=\{f(x)\partial_x\}$, 
$\Dp=\sum f_i(x)\partial_x^i$ 
and the same things with $P$ replaced by $Q$ and $x$ replaced by $y$. 
Denote by $\bDp$ (respectively $\bDq$) the subring of $\Dp$
(respectively $\Dq$) consisting of the operators which annihilate
constant functions.

Note that the space of bivector fields is isomorphic to
$\Vp\otimes_{\R}\Vq$. Let $\D_{P,k}$ be the space of maps
$Vp^{\otimes k}\arrow \bDp$ (which are differential operators
in terms of the coefficients).
\begin{Theorem} There exists a unique sequence $c_k\in \D_{P,k}\otimes
\D_{Q,k}, k=0,1,2,\ldots:$ $c_k=\sum_i a^i_k\otimes b^i_k$, 
$c_0(X,Y)=1\otimes 1$ such that for any 
bivector field $\Psi=\sum_i X_i\wedge Y_i$, $X_i\in \Vp,Y_i\in \Vq$, the formula
\begin{equation}
\begin{array}{ll}
m(\Psi,P,Q,f,g)&=fg+\sum_{k,i_1\ldots i_{k+1}}h^{k+1}L_{X_{i_1}}\{a^n_k(X_{i_2},X_{i_3},\ldots,X_{i_{k+1}})f\}\\ &L_{Y_{i_1}}\{b^n_k(Y_{i_2},Y_{i_3},\ldots,Y_{i_{k+1}})g\}\\
 &=\sum_k m_k(f,g).
\end{array}   \label{posya} 
\end{equation} 
gives a star-product.
\end{Theorem}

{\bf Remark 1.} The ansatz for the formula originates from the following observation. Given a product $m$ from $\Sp$, consider the 'set of zeros' 
of the $m(x,y)$. One can easily show that this set (up to 'biregular
isomorphisms') is an invariant of the star product. Therefore, it is
natural to require that $m(\Psi,x,y)$ would be divisible by $\Psi$. 

To prove this theorem we need some preparation.
Let us pass to the coordinates $x,y$. Then $X_i=\xi_i\partial_x$, $Y_i=\eta_i
\partial_y$, $a^n_k(X_1,\ldots,X_k)f=A^n_k(\xi_1,\ldots,\xi_k,f)$, where
$A_n^k$ are polydifferential operators depending on the derivatives of 
$\xi$ and $f$ with respect to $x$. Similarly,
$b^n_k(Y_1,\ldots,Y_k,g)=B^n_k(\eta_1,\ldots,\eta_k,g)$.  
Our task is to solve the reccurent equation
\be
bm_k=-\frac 12 \sum\limits_{i=1}^{k-1} [m_i,m_{k-i}], \label{rec}
\ee

such that all $m_k$ are of the form in (\ref{posya}). Here $b$ is the Hochschild
differential. Let us specify 
 the meaning of the conditions imposed by (\ref{posya}).  First, note that all $m_i$ 
as well as $[m_i,m_j]$
belong to a subcomplex $K^{\cdot}$ of the Hochshild complex $C^{\cdot}(\bar A,A)$
($A=C^{\infty}(R^2)$ and $\bar A=A/\R$) such that $K^0=K^1=0$; $K^i=$ polydifferential operators $D(f_1,\ldots,f_i),$ such that $D(\phi(y)f_1,\ldots,f_i\psi(x))=\phi(y)\psi(x)D(f_1,\ldots,f_i)$.
That is, we take the cochains that only depend on $\partial_x^m f_1,\partial_y^l f_i$ and  derivatives of $f_2,\ldots,f_{i-1}$.
\begin{Lemma}
The cohomology of $K^{\cdot}$ is generated over $C^{\infty}(R^2)$ by the class of 
$\partial_x\otimes \partial_y$.
\end{Lemma} 

Now, let us specify exactly the space in which all $m_k$ should be. Note that
$$
m_k=\sum\limits_{n,i_1,\ldots,i_k}\xi_{i_ 1}{\partial_x}\circ A^n_{k-1}(\xi_{i_2},\ldots,\xi_{i_k})\otimes\eta_{i_1}{\partial_y}\circ B^n_{k-1}(\eta_{i_2},\ldots,\eta_{i_k}).
$$
Denote by $\ep\in\D_{P,k}$ the space of operators of the form $D=\partial_x\circ D_1$
with $D_1\in\D_{P,k}$. Certainly, $\ep$ depends on a choice of the coordinate $x$ Define $\eq$ in the same fashion.

Recall that $\Psi=\sum_i\xi_i\eta_i\partial_x\wedge\partial_y$. Put
$\phi=\sum_i\xi_i\eta_i$. Then our theorem means exactly that 
\be m_k=\phi K,
\label{met}
\ee
 where $K\in \ep\otimes\eq$. 

Let us investigate how $\ep$ and $\eq$ interact with the Hochshild differential.
Let $\A$ be the subalgebra of functions
depending on the derivatives of $\xi_1,\ldots\xi_k$. Put
$L^i_P=\D_{P,k}^{\otimes_{\A}^i}$ and the similar for $L^i_Q$. Then we have the Hochshild differential
$b:L^i_P\arrow L^{i+1}_P$.
\begin{Lemma}
The sequence
$$
0{\arrow}\ep\stackrel{b}{\arrow}L^2_P\stackrel{b}{\arrow}L^3_P
$$
is exact. This is also true if we replace $P$ by $Q$.
\end{Lemma}
{\em Proof}. Since the one-dimensional Hochshild complex is acyclic
for dimensions bigger than 1,
it suffices to show that $L^1_P=D_{P,k}$ is equal to $\ep\oplus \D^1_{P,k}$, where
$\D^1_{P,k}$ is the subset of operators of order 1 from $\D_{P,k}$. But this 
splitting is given by the Euler-Lagrange operator $E: \D_{P,k}\arrow\D^1_{P,k}$,
$$
E(\sum_{i=1}^N a_i\partial^i_x)=\left(\sum_{i=1}^N(-1)^{i-1}\partial^{i-1}_xa_i\right)\partial_x,
$$ 
since the kernel of $E$ is exactly $\ep$.

{\em Proof of the Theorem 1}. Suppose we have found $m_1,\ldots,m_{k-1}$.
Show that we can solve (\ref{rec}) for $m_k$ so that it is  of the form in (\ref{met}).
Denote by $A$ the right hand side of (\ref{rec}). Note that
$bA=0$ and $A\in K^3$. This means that $A=bS$, where $S\in K^2$.
Note that 
\be
bK^2\in L_P^2\otimes L_Q^1\oplus L_P^1\otimes L_Q^2.\label{dec}
\ee
Let us define a projector $p_P:K^3\arrow L_P^2\otimes L_Q^1$.
For this let us notice that the right hand side can be interpreted
as a differential map from $\C_P\otimes\C_P\otimes\C_Q$ to $C^{\infty}{\R^2}$,
(depending on $\xi_1,\ldots\xi_k,\eta_1,\ldots,\eta_k$), whereas
the left hand side is a differential map from $\C_P\otimes{\C^{\infty}{\R^2}}\otimes
\C_Q$ to $C^{\infty}{\R^2}$. Thus, $ p$ can be defined as a restriction from $C_P\otimes \C^{\infty}{\R^2}
\otimes\C_Q$ to $\C_P\otimes\C_P\otimes\C_Q$. In the same way we can define the
projector $p_Q$ onto the second summand in (\ref{dec}). 
Thus,
$$
\begin{array}{ll} 
A &=(p_P+p_Q)A\\
  &=\sum \phi\partial_x\circ a_p^n(\xi_{i_2},\ldots,\xi_{i_p})
\circ \xi_{j_1}\partial_x\circ a_q^m(\xi_{j_2},\ldots,\xi_{j_q})\\
&\otimes \eta_{j_1}\partial_yb_q^m(\eta_{j_2},\ldots\eta_{j_n})\otimes\partial_y
b_p^n(\eta_{i_1},\ldots,\eta_{i_p})\\
&- \phi\partial_xa_p^n(\xi_{i_2},\ldots\xi_{i_n})\otimes\xi_{j_1}\partial_x
a_q^m(\xi_{j_1},\ldots,\xi_{j_q})\\
&\otimes \partial_y\circ b_p^n(\eta_{i_2},\ldots,\eta_{i_p})
\circ \eta_{j_1}\partial_y\circ b_q^m(\eta_{j_2},\ldots,\eta_{j_q}).
\end{array}
$$
From this we deduce that $A=\phi T$, where $T\in L^2_P\otimes\eq\oplus
\ep\otimes L^2_Q$. Set $T=t_1+t_2$, where $t_1\in L^2_P\otimes\eq$ and
$t_2\in \ep\otimes L^2_Q$. Since $bA=0$, we see that 
$(b\otimes 1) t_1=0$. By Lemma 2, $t_1=(bl_i)\otimes m_i$, where $l_i\in \ep$,
$m_i\in \eq$. Similarly,  $t_2=u_i\otimes bv_i$, where $u_i\in \ep$,
$v_i\in \eq$. Using the closeness of $A$, we have
$-bl_i\otimes bm_i+bu_i\otimes bv_i=0$. Since (by Lemma 2) $b:\ep\arrow L^2$
is an inclusion, we have $c=l_i\otimes m_i=u_i\otimes v_i$ and 
$\phi bc= A$. Thus, we can put $m_k=\phi c$. This proves the existence.
Let us prove the uniqueness. By Lemma 1, the ambiguity in the choice of $m_k$
is of the form $f\partial_x\otimes\partial_y$. This should be of the form
$\phi K$ with $K\in \ep\otimes\eq$, which is impossible, since all non-zero
operators in $\ep,\eq$ have order greater than 1. This proves the theorem. 

The following corollary follows immediately from the uniqueness of the constructed star-product-product. 
\begin{Corollary}
1) The map $(\Psi,P,Q)\mapsto m(\Psi,P,Q)$ is equivariant with respect to Diffeo$(\R^2)$.
2)All operators $a_k^n, b_k^n$ are invariant with respect to the natural action of the group
Diffeo($\R^1$).
\end{Corollary}
It looks very plausible that all the operators $a_k(X_1,X_2,\ldots,X_k)$ are just linear combinations
of $L_{X_{i_1}}\ldots L_{X_{i_k}}$, where $(i_1,\ldots,i_k)$ is a permutation, and the same for $b_k$.
  This is correct
at least for $k\leq 4$.

\begin{Proposition}
The constructed map $p_1: \P\to\Sp$  is a bijection. The inverse map $p_2: \Sp\to \P$
is a well defined map.
\end{Proposition} 
{\em Proof.} Consider the map $p_3:\Sp\to\P:\ m\mapsto\displaystyle\frac{1}{h} (m(x,y)-m(y,x))\partial_x\wedge\partial_y$.
This is a well defined
map, and it is not  hard to see that $p_3$ is injective. Note that  $p_3\circ p_1(\psi)=\psi+O(h)$, hence,
$p_3\circ p_1$ is invertible, and we can put $p_2=(p_3\circ p_1)^{-1}\circ p_3$.   \medskip

Let us check the properties (\ref{equiv}). Let $D_t=\exp(tX)$ be a one-parameter local Lie group of diffeomorphisms corresponding to a vector field $X$. Then $D_t$ acts naturally on $\S$. Using Corollary from the Theorem 1
and Proposition, we can write
\be 
p_1(D_t\Psi)=D_t(m(\psi,D_t^{-1}P,D_t^{-1}Q))\equiv \mp(m(\psi,D_t^{-1}P,D_t^{-1}Q))=
m(\chi_t,P,Q)\label{pp},
\ee
where $\chi_t=p_2 m(\psi,D_t^{-1}P,D_t^{-1}Q)$. Put $\delta_X(\Psi)=\displaystyle\frac{d}{dt}\chi_t|_{t=0}$.
This is a well defined map, linear in $X$. Furthermore, $\dX=O(h)$. It is enough to prove the following.
\begin{Theorem}. There exists  a linear differential operator $A:\Vect[[h]]\to\Vect[[h]]$, depending
on $\Psi\in\P$, such that
$$
\dX(\Psi)=L_{A(X)}\Psi.
$$
\end{Theorem} 
Let us explain why it is enough. First, disregarding, if needed, some terms, we can make $A(X)$
to be $O(h)$. Then (\ref{pp}) can be rewritten as $p_1(\Psi+tL_X\Psi)\equiv p_1(\Psi+tL_{A(x)}\Psi)+o(t)$,
$p_1(\Psi+tL_{X-A(X)}\Psi)\equiv p_1(\Psi)+o(t)$. Since $A(X)=O(h)$, the equation
$Y=X-A(X)$ is solvable for all $Y$ and we deduce that $p_1(D_t\Psi)\equiv p_1(\Psi)$. Since $p_1$ is
automatically equivariant with respect to the reflection $(x,y)\to (-x,y)$, this would mean that 
$p_1$ is equivariant with respect to the whole group of diffeomorphisms. Thus, $p_1(x)\equiv p_1(y)$.
Suppose that $m_1\equiv m_2$. Then $p_2(m_{1,2})=p_2(\mu_{1,2})$, where $\mu_{1,2}=\mp(m_{1,2})$. Then
$\mu_2(f,g)=UD\mu_1((UD)^{-1}f,(UD)^{-1}g)$, where $D$ is a formal diffeomorphism
(that is, an element of $\rm Diffeo\stackrel{.}{+}\exp h{\rm Vect}[[h]]$ ) and $U$ is a differential operator $U=1+hV$, 
satisfying the condition 1 of Proposition 1. Thus, $p_2(\mu_2)=p_2(D^{-1}\mu_1)$. If $D^{-1}=\exp(tX)$, 
and $p_2(\mu_1)=\Psi$, then $p_2\mu_2=\chi_t\equiv \Psi$. We only need to check the equivariance with respect
to the reflection $(x,y)\mapsto (-x,y)$, which immediately[C[ follows from the Corollary 1.
Thus, we only need to prove  Theorem 2. We need to make a reduction.
Denote by $\V$ the space of all linear over $\R[[h]]$ differential operators $\Vect[[h]]\to\Vect[[h]]$,
depending on $\Psi\in \P$. Let $\Psi=\phi\partial_x\wedge\partial_y$.
Denote $\Vr=\overline{\V[\phi^{-1}]}$, where the bar means the completion in 
the $h$-adic topology.
\begin{Proposition} It is enough to find $A$ in $\Vr$.
\end{Proposition}
{\em Proof.} Suppose we have found such an $A$. Let $l$ be the least degree in
$h$, where A has singularity. Let $N$
be the least positive integer such that $B=\phi^NA$ does not
have singularities up to $h^{l+1}$. Further we will write
$a\equiv b$ if $\phi^{N-1}(a-b)$ does not have singularities up to $h^{l+1}$.
Put $B=U\partial_x+V\partial_y$. Then 
$$
0=L_A\Psi\equiv (1+n)\frac{U\phi_x+V\phi_y}{\phi^N}(\partial_x\wedge\partial_y).
$$
Hence, 
$$
B\equiv \frac{W(\phi_x\partial_y-\phi_y\partial_x)}{\phi^n}\equiv
\frac{1}{N}\{\frac{W}{\phi^N},\Psi\},
$$
where $\{,\}$ is the Schouten bracket. Put
$A_1=A-\frac{1}{N}\{\frac{W}{\phi^N},\Psi\}$. Then $L_{A_1(X)}\Psi=L_{A(X)}\Psi$
and $A_1$ has no singularities of order $N$ up to $h^{l+1}$. Iterating this
procedure, we will get rid of all the singularities.

Let us make the following reductions.
Similarly to the differenrial operators on vector fields put
$\F=\overline{\Fun(\P)[\phi^{-1}]}$, where $\Fun(\P)$ is the space of functions on $\P$. For $F\in \F$ put
$\dX F(\Psi)=\displaystyle\frac{d}{dt}F(\Psi+t\dX\Psi)|_{t=0}$. 
Put $\omega=\displaystyle\frac{1}{\phi}dx\wedge dy$. 
Then we only need to prove that 
$$
\dX\omega=d\theta(\Psi,X).
$$ for some 1-form $\theta$.
It is clear that it is enough  to prove this for some form 
\be
\Omega=\omega+d\alpha. \label{exac}
\ee Let us find a suitable form $\Omega$.
Also, we can assume that out vector field $X$ is tangent to $Q$ (since any vector field is a sum of a vector field 
tangent to $P$ and a vector field tangent to $Q$).
\begin{Proposition}
\begin{enumerate} Let $z$ be some function on $R^2$ such that $(x,z)$ form a nondegenerate coordinate system.
Put $\Psi=\phi\partial_x\wedge\partial_z$.
\item[1)]For any $m=m(\Psi,P,Q)$ the differential operator $\displaystyle\frac{1}{h}\ad x$ can be represented as 
$\phi\partial_z(1+ S(m)\circ\partial_z)$ for some well defined differential operator $S=S(\Psi,P,Q)=O(h)$.
\item[2)]There exists a unique $f=f(x,z,P,Q,\Psi)\in \F $ such that 
\be
\phi(1+\partial_z\circ S) f =1.\label{it2}
\ee
\item[3)] Put $\Omega=\Omega(x,z,P,Q,\Psi)=fdx\wedge dz$. Then  (\ref{exac}) holds.
\end{enumerate}
\end{Proposition}
{\em Proof.} 1)This immediately follows from (\ref{posya}). The second statement holds because $S=O(h)$. The last statement is true
because (\ref{it2}) can be rewritten as 
\be f=1/\phi+\partial_z\tau\label{part}
\ee for some $\tau$.

{\bf Remark.} If we had an antiderivative F of $f$, such that $F_y=f$, then  it would be $[x,F]=h$ and
$\Omega=dx\wedge dF$. Thus, $\Omega$ is nothing else but the Berezin curvature
\cite{Ber}.
 Also, it is not hard to
prove that our construction does not depend on a choice of $z$. 

Now we are ready to prove the invariance. Formula (\ref{pp})
can be rewritten as
$$
m_{P,Q_t}m(\chi_t,P,Q)=m(\Psi,P,Q_t),
$$
where $Q_t=D_t^{-1}Q$.
Let $A_t\circ\partial_z=1/h\ad_{m(\chi_t,P,Q)}x$,   
$B_t\circ\partial_z=1/h\ad_{m(\Psi,P,Q_t)}x$.
Recall that $m_{P,Q_t}$ is a conjugation with respect to some operator $U_t=1+hV_t$, and $U_t(x^n)=x^n $(see (\ref{polar})).
Therefore, $U_t=1+hW_t\partial_z$. Also, $B_t\partial_z=U_tA(t)\partial_zU^{-1}_t$ and $A_tf(\chi_t,P,Q)=B_t(f(\Psi,x,z,P,Q_t))=1$. One can check that $f(\Psi,x,z,P,Q_t)=f(\chi_t,P,Q)+\partial_z\circ
hW_tf(\chi_t,P,Q)$. 
Using (\ref{part}), we  immediately get $\omega=\chi_t^{-1}+\partial_z\tau_t dx\wedge dz$, where
$\chi_t^{-1}=(\chi_t,dx\wedge dy)^{-1}dx\wedge dy$. Therefore, $\dX \omega=d\alpha$, where
$$
\alpha=\frac{d}{dt}d(\tau_tdx)|_{t=0}.
$$
 Which proves theorem 2.

\end{document}